\documentclass[11pt]{amsart}
\usepackage{amsmath,amssymb,mathrsfs}
\newtheorem{theorem}{Theorem}
\newtheorem{proposition}[theorem]{Proposition}
\newtheorem{corollary}[theorem]{Corollary}
\newtheorem{lemma}[theorem]{Lemma}

\begin{document}

\title{Rigidity of area-minimizing two-spheres in three-manifolds}
\author{H. Bray, S. Brendle, and A. Neves}
\address{Department of Mathematics \\ Duke University \\ Durham, NC 27708}
\address{Department of Mathematics \\ Stanford University \\ 450 Serra Mall, Bldg 380 \\ Stanford, CA 94305} 
\address{Imperial College \\ Huxley Building \\ 180 Queen's Gate \\ London SW7 2RH \\ United Kingdom}

\maketitle

\section{Introduction}

A classical result in differential geometry due to Toponogov \cite{Topogonov} states that every simple closed geodesic $\gamma$ on a two-dimensional surface $(\Sigma,g)$ satisfies
\[\text{\rm length}(\gamma)^2 \, \inf_\Sigma K \leq 4\pi^2,\] 
where $K$ denotes the Gaussian curvature of $\Sigma$. Moreover, equality holds if and only if $(\Sigma,g)$ is isometric to the standard sphere $S^2$ up to scaling (see also \cite{Hang-Wang} for a different proof).

We next consider a three-manifold $(M,g)$ with positive scalar curvature. By a theorem of Schoen and Yau \cite{Schoen-Yau}, any area-minimizing surface in $M$ is homeomorphic to either $S^2$ or $\mathbb{RP}^2$. The case of area-minimizing projective planes was studied in \cite{Bray-Brendle-Eichmair-Neves}. In particular, if $\Sigma$ is an embedded $\mathbb{RP}^2$ of minimal area, then the area of $\Sigma$ can be estimated from above by 
\[\text{\rm area}(\Sigma,g) \, \inf_M R \leq 12\pi\] 
(cf. \cite{Bray-Brendle-Eichmair-Neves}, Theorem 1). Moreover, equality holds if and only if $(M,g)$ is isometric to $\mathbb{RP}^3$ up to scaling.

In this paper, we consider the case of area-minimizing two-spheres. We shall assume throughout that $(M,g)$ is a compact three-manifold with $\pi_2(M) \neq 0$. We denote by $\mathscr{F}$ the set of all smooth maps $f: S^2 \to M$ which represent a non-trivial element of $\pi_2(M)$. We define
\begin{equation} 
\label{definition.of.A}
\mathscr{A}(M,g) = \inf \{\text{\rm area}(S^2,f^* g): f \in \mathscr{F}\}. 
\end{equation}
We now state the main result of this paper:

\begin{theorem}
\label{main.theorem}
We have
\begin{equation}
\label{upper.bound.for.area}
\mathscr{A}(M,g) \, \inf_M R \leq 8\pi, 
\end{equation}
where $R$ denotes the scalar curvature of $(M,g)$. Moreover, if equality holds, then the universal cover of $(M,g)$ is isometric to the standard cylinder $S^2 \times \mathbb{R}$ up to scaling.
\end{theorem}

The inequality \eqref{upper.bound.for.area} follows directly from the formula for the second variation of area. We now describe the proof of the rigidity statement. By scaling, we may assume that $\mathscr{A}(M,g) = 4\pi$ and $\inf_M R = 2$. It follows from results of Meeks and Yau \cite{Meeks-Yau} that the infimum in \eqref{definition.of.A} is attained by a smooth immersion $f \in \mathscr{F}$ (see also \cite{Hass-Scott}). Using the implicit function theorem, we construct a one-parameter family of immersed two-spheres with constant mean curvature. Using the formula for the second variation of area, we are able to show that these surfaces all have area $\mathscr{A}(M,g) = 4\pi$. Consequently, these spheres are all round and totally geodesic. This allows us to construct a local isometry from the cylinder $S^2 \times \mathbb{R}$ into $M$. The use of constant mean curvature surfaces is motivated in part by work of Bray \cite{Bray-thesis} and Huisken and Yau \cite{Huisken-Yau} (see also \cite{Bray-Penrose} and \cite{Huisken-Ilmanen}).

We note that Cai and Galloway \cite{Cai-Galloway} have obtained a similar rigidity theorem for minimal tori in three-manifolds of nonnegative scalar curvature. The proof in \cite{Cai-Galloway} uses a different argument based on a deformation of the metric to strictly positive scalar curvature. The arguments in this paper can be adapted to give an alternative proof of Theorem 1 in \cite{Cai-Galloway}. See also \cite{Galloway} for related work in this direction.

The authors would like to thank Fernando Marques for discussions, and Willie Wong for comments on an earlier version of this paper.

\section{Proof of (\ref{upper.bound.for.area})}

Let us consider a smooth immersion $f: S^2 \to M$. Since $f$ has trivial normal bundle, there exists a globally defined unit normal field $\nu$. In other words, for each point $x \in S^2$, $\nu(x) \in T_{f(x)} M$ is a unit vector which is orthogonal to the image of $df_x: T_x S^2 \to T_{f(x)} M$. The following result is a consequence of the Gauss-Bonnet theorem:

\begin{proposition} 
\label{gauss.bonnet}
For any immersion $f: S^2 \to M$, we have 
\[\int_{S^2} (R - 2 \, \text{\rm Ric}(\nu,\nu) - |I\!I|^2) \, d\mu_{f^* g} \leq 8\pi,\]
where $I\!I$ denotes the second fundamental form of $f$.
\end{proposition} 

\textbf{Proof.} 
By the Gauss equations, we have 
\[R - 2 \, \text{\rm Ric}(\nu,\nu) - |I\!I|^2 = 2K - H^2,\] 
where $H$ and $K$ denote the mean curvature and the Gaussian curvature, respectively. This implies 
\[\int_{S^2} (R - 2 \, \text{\rm Ric}(\nu,\nu) - |I\!I|^2) \, d\mu_{f^* g} \leq 2 \int_{S^2} K \, d_{f^* g} = 8\pi\] 
by the Gauss-Bonnet theorem. \\

We next consider a map $f \in \mathscr{F}$ which attains the infimum in (\ref{definition.of.A}). The existence of a minimizer is guaranteed by the following result:

\begin{proposition}
\label{existence.of.minimizer}
There exists a smooth map $f \in \mathscr{F}$ such that $\text{\rm area}(S^2,f^* g) = \mathscr{A}(M,g)$. Moreover, the map $f$ is an immersion.
\end{proposition}

Proposition \ref{existence.of.minimizer} is an immediate consequence of Theorem 7 in \cite{Meeks-Yau} (see also \cite{Hass-Scott}, Theorem 4.2). In fact, Meeks and Yau show that either $f$ is an embedding or a two-to-one covering map whose image is an embedded $\mathbb{RP}^2$. We will not use this stronger statement here.

Let $f \in \mathscr{F}$ be a smooth immersion with $\text{\rm area}(S^2,f^* g) = \mathscr{A}(M,g)$. Using the formula for the second variation of area, we obtain
\[\int_{S^2} (\text{\rm Ric}(\nu,\nu) + |I\!I|^2) \, u^2 \, d\mu_{f^* g} \leq \int_{S^2} |\nabla u|_{f^* g}^2 \, d\mu_{f^* g}\]
for every smooth function $u: S^2 \to \mathbb{R}$. Choosing $u = 1$ gives 
\[\int_{S^2} (\text{\rm Ric}(\nu,\nu) + |I\!I|^2) \, d\mu_{f^* g} \leq 0.\] 
Using Proposition \ref{gauss.bonnet}, we obtain 
\begin{align}
\label{estimate}
\text{\rm area}(S^2,f^* g) \, \inf_M R 
&\leq \int_{S^2} (R + |I\!I|^2) \, d\mu_{f^* g} \notag \\
&\leq 8\pi + 2 \int_{S^2} (\text{\rm Ric}(\nu,\nu) + |I\!I|^2) \, d\mu_{f^* g} \\
&\leq 8\pi. \notag
\end{align}
This completes the proof of (\ref{upper.bound.for.area}).

\section{The case of equality}

In this section, we analyze the case of equality. Suppose that $$\mathscr{A}(M,g) \, \inf_M R = 8\pi.$$ After rescaling the metric if necessary, we may assume that $\mathscr{A}(M,g) = 4\pi$ and $\inf_M R = 2$. By Proposition \ref{existence.of.minimizer}, we can find a smooth immersion $f \in \mathscr{F}$ such that $\text{\rm area}(S^2,f^* g) = 4\pi$. 

\begin{proposition}
\label{round}
Let $f \in \mathscr{F}$ be a smooth immersion such that $\text{\rm area}(S^2,f^* g) = 4\pi$. Then the surface $\Sigma = f(S^2)$ is totally geodesic. Moreover, we have $R = 2$ and $\text{\rm Ric}(\nu,\nu) = 0$ at each point on $\Sigma$.
\end{proposition}

\textbf{Proof.}
By assumption, we have $\text{\rm area}(S^2,f^* g) = 4\pi$ and $\inf_M R = 2$. Consequently, the inequalities in (\ref{estimate}) are all equalities. In particular, we have 
\begin{equation}
\label{a}
\int_{S^2} (R + |I\!I|^2) \, d\mu_{f^* g} = 8\pi
\end{equation} 
and 
\begin{equation} 
\label{b} 
\int_{S^2} (\text{\rm Ric}(\nu,\nu) + |I\!I|^2) \, d\mu_{f^* g} = 0. 
\end{equation}
It follows from (\ref{b}) that the constant functions lie in the nullspace of the Jacobi operator $L = \Delta_{f^* g} + \text{\rm Ric}(\nu,\nu) + |I\!I|^2$. This implies
\[\text{\rm Ric}(\nu,\nu) + |I\!I|^2 = 0\]
at each point on $\Sigma$. Moreover, since $\text{\rm area}(S^2,f^* g) = 4\pi$ and $\inf_M R = 2$, the identity (\ref{a}) implies that $R = 2$ and $|I\!I|^2 = 0$ at each point on $\Sigma$. This completes the proof. \\

\begin{proposition}
\label{existence.of.cmc.foliation}
Let $f \in \mathscr{F}$ be a smooth immersion such that $\text{\rm area}(S^2,f^* g) = 4\pi$. Then there exists a positive real number $\delta_1$ and a smooth map $w: S^2 \times (-\delta_1,\delta_1) \to \mathbb{R}$ with the following properties:
\begin{itemize}
\item For each point $x \in S^2$, we have $w(x,0) = 0$ and $\frac{\partial}{\partial t} w(x,t) \big |_{t=0} = 1$.
\item For each $t \in (-\delta_1,\delta_1)$, we have $\int_{S^2} (w(\cdot,t) - t) \, d\mu_{f^* g} = 0$.
\item For each $t \in (-\delta_1,\delta_1)$, the surface 
\[\Sigma_t = \{\exp_{f(x)}(w(x,t) \, \nu(x)): x \in S^2\}\] 
has constant mean curvature.
\end{itemize}
\end{proposition}

\textbf{Proof.} 
The Jacobi operator associated with the minimal immersion $f: S^2 \to M$ is given by $L = \Delta_{f^* g} + \text{\rm Ric}(\nu,\nu) + |I\!I|^2$. Using Proposition \ref{round}, we conclude that $L = \Delta_{f^* g}$. Hence, the assertion follows from the implicit function theorem. \\

For each $t \in (-\delta_1,\delta_1)$, we define a map $f_t: S^2 \to M$ by $f_t(x) = \exp_{f(x)}(w(x,t) \, \nu(x))$. Clearly, $f_0(x) = f(x)$ for all $x \in S^2$. To fix notation, we denote by $\nu_t(x) \in T_{f_t(x)} M$ the unit normal vector to the surface $\Sigma_t = f_t(S^2)$ at the point $f_t(x)$. We assume that $\nu_t$ depends smoothly on $x$ and $t$, and $\nu_0(x) = \nu(x)$ for all $x \in S^2$. Moreover, we denote by $I\!I_t$ the second fundamental form of $f_t$.

\begin{lemma}
\label{prep.1}
There exists a positive real number $\delta_2 < \delta_1$ with the following property: if $t \in (-\delta_2,\delta_2)$ and $u: S^2 \to \mathbb{R}$ is a smooth function satisfying $\int_{S^2} u \, d\mu_{f_t^* g} = 0$, then 
\[\int_{S^2} |\nabla u|_{f_t^* g}^2 \, d\mu_{f_t^* g} - \int_{S^2} (\text{\rm Ric}(\nu_t,\nu_t) + |I\!I_t|^2) \, u^2 \, d\mu_{f_t^* g} \geq 0.\]
\end{lemma}

\textbf{Proof.} 
We can find a uniform constant $c > 0$ such that 
\[\int_{S^2} |\nabla u|_{f_t^* g}^2 \, d\mu_{f_t^* g} \geq c \int_{S^2} u^2 \, d\mu_{f_t^* g}\] 
for each $t \in (-\delta_1,\delta_1)$ and every smooth function $u: S^2 \to \mathbb{R}$ satisfying $\int_{S^2} u \, d\mu_{f_t^* g} = 0$. Moreover, it follows from Proposition \ref{round} that 
\[\sup_{S^2} (\text{\rm Ric}(\nu_t,\nu_t) + |I\!I_t|^2) \to 0\] 
as $t \to 0$. Putting these facts together, the assertion follows. \\

\begin{lemma}
\label{prep.2}
For each $t \in (-\delta_1,\delta_1)$, we have 
\[\int_{S^2} (\text{\rm Ric}(\nu_t,\nu_t) + |I\!I_t|^2) \, d\mu_{f_t^* g} \geq 0.\]
\end{lemma} 

\textbf{Proof.} 
Since $f$ minimizes area in its homotopy class, we have 
\[\text{\rm area}(S^2,f_t^* g) \geq \text{\rm area}(S^2,f^* g) = 4\pi.\] 
Moreover, we have $\inf_M R = 2$. Applying Proposition \ref{gauss.bonnet} to the map $f_t: S^2 \to M$, we obtain 
\begin{align*}
8\pi
&\leq \text{\rm area}(S^2,f_t^* g) \, \inf_M R \\
&\leq \int_{S^2} (R + |I\!I_t|^2) \, d\mu_{f_t^* g} \\
&\leq 8\pi + 2 \int_{S^2} (\text{\rm Ric}(\nu_t,\nu_t) + |I\!I_t|^2) \, d\mu_{f_t^* g}.
\end{align*} From this, the assertion follows. \\

By assumption, the surface $\Sigma_t$ has constant mean curvature. The mean curvature vector of $\Sigma_t$ can be written in the form $-H(t) \, \nu_t$, where $H(t)$ is a smooth function of $t$. For each $t \in (-\delta_1,\delta_1)$, the lapse function $\rho_t: S^2 \to \mathbb{R}$ is defined by 
\begin{equation} 
\label{lapse.function}
\rho_t(x) = \Big \langle \nu_t(x),\frac{\partial}{\partial t} f_t(x) \Big \rangle. 
\end{equation}
Clearly, $\rho_0(x) = 1$ for all $x \in S^2$. By continuity, we can find a positive real number $\delta_3 < \delta_2$ such that $\rho_t(x) > 0$ for all $x \in S^2$ and all $t \in (-\delta_3,\delta_3)$. The lapse function $\rho_t: S^2 \to \mathbb{R}$ satisfies the Jacobi equation 
\begin{equation} 
\label{Jacobi.equation}
\Delta_{f_t^* g} \rho_t + (\text{\rm Ric}(\nu_t,\nu_t) + |I\!I_t|^2) \, \rho_t = -H'(t) 
\end{equation}
(cf. \cite{Huisken-Ilmanen}, equation (1.2)). 

\begin{proposition}
\label{area}
We have $\text{\rm area}(S^2,f_t^* g) = 4\pi$ for all $t \in (-\delta_3,\delta_3)$.
\end{proposition}

\textbf{Proof.} 
Let $\overline{\rho}_t$ denote the mean value of the lapse function $\rho_t: S^2 \to \mathbb{R}$ with respect to the induced metric $f_t^* g$; that is, 
\[\overline{\rho}_t = \frac{1}{\text{\rm area}(S^2,f_t^* g)} \int_{S^2} \rho_t \, d\mu_{f_t^* g}.\] 
It follows from Lemma \ref{prep.1} that 
\[\int_{S^2} |\nabla \rho_t|_{f_t^* g}^2 \, d\mu_{f_t^* g} - \int_{S^2} (\text{\rm Ric}(\nu_t,\nu_t) + |I\!I_t|^2) \, (\rho_t - \overline{\rho}_t)^2 \, d\mu_{f_t^* g} \geq 0\] 
for all $t \in (-\delta_2,\delta_2)$. Moreover, Lemma \ref{prep.2} implies that 
\[\overline{\rho}_t^2 \int_{S^2} (\text{\rm Ric}(\nu_t,\nu_t) + |I\!I_t|^2) \, d\mu \geq 0\] 
for all $t \in (-\delta_1,\delta_1)$. Adding both inequalities yields 
\begin{equation} 
\label{derivative.of.H.1}
\int_{S^2} |\nabla \rho_t|_{f_t^* g}^2 \, d\mu_{f_t^* g} + \int_{S^2} (\text{\rm Ric}(\nu_t,\nu_t) + |I\!I_t|^2) \, \rho_t \, (2 \, \overline{\rho}_t - \rho_t) \, d\mu_{f_t^* g} \geq 0 
\end{equation}
for all $t \in (-\delta_2,\delta_2)$.

In the next step, we multiply the equation (\ref{Jacobi.equation}) by $2 \, \overline{\rho}_t - \rho_t$ and integrate. This gives 
\begin{align} 
\label{derivative.of.H.2}
&\int_{S^2} |\nabla \rho_t|_{f_t^* g}^2 \, d\mu_{f_t^* g} + \int_{S^2} (\text{\rm Ric}(\nu_t,\nu_t) + |I\!I_t|^2) \, \rho_t \, (2 \, \overline{\rho}_t - \rho_t) \, d\mu_{f_t^* g} \notag \\ 
&= -H'(t) \int_{S^2} (2 \, \overline{\rho}_t - \rho_t) \, d\mu_{f_t^* g} = -H'(t) \int_{S^2} \rho_t \, d\mu_{f_t^* g}. 
\end{align}
Putting these facts together, we obtain 
\begin{equation} 
\label{derivative.of.H.3}
H'(t) \int_{S^2} \rho_t \, d\mu_{f_t^* g} \leq 0 
\end{equation}
for each $t \in (-\delta_2,\delta_2)$. Therefore, we have $H'(t) \leq 0$ for all $t \in (-\delta_3,\delta_3)$. Since $H(0) = 0$, it follows that $H(t) \geq 0$ for all $t \in (-\delta_3,0]$ and $H(t) \leq 0$ for all $t \in [0,\delta_3)$. Using the identity 
\[\frac{d}{dt} \text{\rm area}(S^2,f_t^* g) = \int_{S^2} \Big \langle H(t) \, \nu_t,\frac{\partial}{\partial t} f_t \Big \rangle \, d\mu_{f_t^* g} = H(t) \int_{S^2} \rho_t \, d\mu_{f_t^* g},\] 
we obtain 
\[\text{\rm area}(S^2,f_t^* g) \leq \text{\rm area}(S^2,f^* g) = 4\pi\] 
for all $t \in (-\delta_3,\delta_3)$. Since $f$ minimizes area in its homotopy class, we conclude that $\text{\rm area}(S^2,f_t^* g) = 4\pi$ for all $t \in (-\delta_3,\delta_3)$. \\

\begin{proposition}
\label{rigidity}
For each $t \in (-\delta_3,\delta_3)$, the surface $\Sigma_t$ is totally geodesic, and we have $R = 2$ and $\text{\rm Ric}(\nu_t,\nu_t) = 0$ at each point on $\Sigma_t$. Moreover, the lapse function $\rho_t: S^2 \to \mathbb{R}$ is constant. Moreover
\end{proposition}

\textbf{Proof.}
By Proposition \ref{area}, we have $\text{\rm area}(S^2,f_t^* g) = 4\pi$. Hence, it follows from Proposition \ref{round} that $\Sigma_t$ is totally geodesic, and $R = 2$ and $\text{\rm Ric}(\nu_t,\nu_t) = 0$ at each point on $\Sigma_t$. Substituting this into (\ref{Jacobi.equation}), we obtain $\Delta_{f_t^* g} \rho_t = 0$. Therefore, the function $\rho_t: S^2 \to \mathbb{R}$ is constant, as claimed. \\

\begin{corollary}
\label{local.splitting} 
The normal vector field $\nu_t$ is a parallel vector field near $\Sigma$. In particular, each point on $\Sigma$ has a neighborhood which is isometric to a Riemannian product.
\end{corollary}

\textbf{Proof.} 
By Proposition \ref{rigidity}, the lapse function $\rho_t: S^2 \to \mathbb{R}$ is constant. This implies 
\begin{align*} 
&\Big \langle D_{\frac{\partial f_t}{\partial x_i}} \nu_t,\frac{\partial f_t}{\partial t} \Big \rangle - \Big \langle D_{\frac{\partial f_t}{\partial t}} \nu_t,\frac{\partial f_t}{\partial x_i} \Big \rangle \\ 
&= \frac{\partial}{\partial x_i} \Big \langle \nu_t,\frac{\partial f_t}{\partial t} \Big \rangle - \frac{\partial}{\partial t} \Big \langle \nu_t,\frac{\partial f_t}{\partial x_i} \Big \rangle = \frac{\partial}{\partial x_i} \rho_t(x) = 0 
\end{align*} 
for each point $x \in S^2$. Moreover, we have 
\begin{align*} 
&\Big \langle D_{\frac{\partial f_t}{\partial x_i}} \nu_t,\frac{\partial f_t}{\partial x_j} \Big \rangle - \Big \langle D_{\frac{\partial f_t}{\partial x_j}} \nu_t,\frac{\partial f_t}{\partial x_i} \Big \rangle \\ 
&= \frac{\partial}{\partial x_i} \Big \langle \nu_t,\frac{\partial f_t}{\partial x_j} \Big \rangle - \frac{\partial}{\partial x_j} \Big \langle \nu_t,\frac{\partial f_t}{\partial x_i} \Big \rangle = 0 
\end{align*}
for all $x \in S^2$. Putting these facts together, we obtain 
\[\Big \langle D_{\frac{\partial f_t}{\partial x_i}} \nu_t,V \Big \rangle - \Big \langle D_V \nu_t,\frac{\partial f_t}{\partial x_i} \Big \rangle = 0\] 
for each point $x \in S^2$ and all vectors $V \in T_{f(x)} M$. In particular, we have 
\[\Big \langle D_{\frac{\partial f_t}{\partial x_i}} \nu_t,\nu_t \Big \rangle - \Big \langle D_{\nu_t} \nu_t,\frac{\partial f_t}{\partial x_i} \Big \rangle = 0\] 
for each point $x \in S^2$. Since the vector field $\nu_t$ has unit length, we conclude that $D_{\nu_t} \nu_t = 0$. On the other hand, it follows from Proposition \ref{rigidity} that the surfaces $\Sigma_t$ are totally geodesic. This implies $D_{\frac{\partial f_t}{\partial x_i}} \nu_t = 0$ for each point $x \in S^2$. Thus, we conclude that the normal vector field $\nu_t$ is parallel. This completes the proof of Corollary \ref{local.splitting}. \\

We now consider the product $S^2 \times \mathbb{R}$, where $S^2$ is equipped with the induced metric $f^* g$. We define a map $\Phi: S^2 \times \mathbb{R} \to M$ by $\Phi(x,t) = \exp_{f(x)}(t \, \nu(x))$. It follows from Corollary \ref{local.splitting} that the restriction $\Phi|_{S^2 \times (-\delta,\delta)}$ is a local isometry if $\delta>0$ is sufficiently small.

\begin{proposition}
\label{isometry}
The map $\Phi: S^2 \times \mathbb{R} \to M$ is a local isometry.
\end{proposition} 

\textbf{Proof.}
We first show that $\Phi|_{S^2 \times [0,\infty)}$ is a local isometry. Suppose this is false. Let $\tau$ be the largest positive real number with the property that $\Phi|_{S^2 \times [0,\tau]}$ is a local isometry. We now define a map $\tilde{f}: S^2 \to M$ by $\tilde{f}(x) = \Phi(x,\tau)$. Clearly, $\tilde{f}$ is homotopic to $f$; consequently, $\tilde{f}$ represents a non-trivial element of $\pi_2(M)$. Moreover, we have $\text{\rm area}(S^2,\tilde{f}^* g) = \text{\rm area}(S^2,f^* g) = 4\pi$. Therefore, $\tilde{f}$ has minimal area among all maps in $\mathscr{F}$. By Corollary \ref{local.splitting}, each point on the surface $\tilde{\Sigma} = \tilde{f}(S^2)$ has a neighborhood which is isometric to a Riemannian product. Hence, if $\delta>0$ is sufficiently small, then the map $\Phi|_{S^2 \times [0,\tau+\delta)}$ is a local isometry. This contradicts the maximality of $\tau$.

Therefore, the restriction $\Phi|_{S^2 \times [0,\infty)}$ is a local isometry. An analogous argument shows that $\Phi|_{S^2 \times (-\infty,0]}$ is a local isometry. This completes the proof of Proposition \ref{isometry}. \\

Since $\Phi: S^2 \times \mathbb{R} \to M$ is a local isometry, it follows that $\Phi$ is a covering map (cf. \cite{Cheeger-Ebin}, Section 1.11). Consequently, the universal cover of $(M,g)$ is isometric to $S^2 \times \mathbb{R}$, equipped with the standard metric. This completes the proof of Theorem \ref{main.theorem}.


\begin{thebibliography}{99}
\bibitem{Bray-thesis}
H.~Bray, \textit{The Penrose inequality in general relativity and volume comparison theorems involving scalar curvature,} PhD thesis, Stanford University (1997)

\bibitem{Bray-Penrose}
H.~Bray, \textit{Proof of the Riemannian Penrose inequality using the positive mass theorem,} J. Diff. Geom. 59, 177--267 (2001)

\bibitem{Bray-Brendle-Eichmair-Neves}
H.~Bray, S.~Brendle, M.~Eichmair, and A.~Neves, \textit{Area-minimizing projective planes in three-manifolds,} Comm. Pure Appl. Math. 63, 1237--1247 (2010)

\bibitem{Cai-Galloway}
M.~Cai and G.~Galloway, \textit{Rigidity of area minimizing tori in $3$-manifolds of nonnegative scalar curvature,} Comm. Anal. Geom. 8, 565--573 (2000)

\bibitem{Cheeger-Ebin}
J.~Cheeger and D.~Ebin, \textit{Comparison theorems in Riemannian geometry,} American Mathematical Society (2008)

\bibitem{Galloway}
G.~Galloway, \textit{Rigidity of marginally trapped surfaces and the topology of black holes,} Comm. Anal. Geom. 16, 217--229 (2008)

\bibitem{Hang-Wang}
F.~Hang and X.~Wang, \textit{Rigidity theorems for compact manifolds with boundary and positive Ricci curvature,} J. Geom. Anal. 19, 628--642 (2009)

\bibitem{Hass-Scott}
J.~Hass and P.~Scott, \textit{The existence of least area surfaces in $3$-manifolds,} Trans. Amer. Math. Soc. 310, 87--114 (1988)

\bibitem{Huisken-Ilmanen}
G.~Huisken and T.~Ilmanen, \textit{The inverse mean curvature flow and the Riemannian Penrose inequality,} J. Diff. Geom. 59, 353--437 (2001)

\bibitem{Huisken-Yau}
G.~Huisken and S.T.~Yau, \textit{Definition of center of mass for isolated physical systems and unique foliations by stable spheres with constant mean curvature,} Invent. Math. 124, 281--311 (1996)

\bibitem{Meeks-Yau}
W.~Meeks and S.T.~Yau, \textit{Topology of three-dimensional manifolds and the embedding problems in minimal surface theory,} Ann. of Math. 112, 441--484 (1980)

\bibitem{Schoen-Yau}
R.~Schoen and S.T.~Yau, \textit{Existence of incompressible minimal surfaces and the topology of three-dimensional manifolds with nonnegative scalar curvature,} Ann. of Math.  110, 127--142 (1979) 

\bibitem{Topogonov}
V.~Toponogov, \textit{Evaluation of the length of a closed geodesic on a convex surface,} Dokl. Akad. Nauk. SSSR 124, 282--284 (1959)
\end{thebibliography}
\end{document}